\documentclass[12pt]{article}
\usepackage{amsfonts}

\newtheorem{thm}{Theorem}[section]

\newtheorem{prop}[thm]{Proposition}
\newtheorem{conjecture}{Conjecture}

\newenvironment{remark}{\par\medskip\noindent{\bf Remark.\ }}{\par\smallskip}
\newcommand{\proof
}{\par\medskip\noindent {\bf Proof.\ \ }}

\newcommand{\be}{\begin{equation}}
\newcommand{\ee}{\end{equation}}
\newcommand{\openbox}{\leavevmode
  \hbox to8pt{\hfil\vrule\vbox to6pt{\hrule width6pt\vfil\hrule}\vrule}}

\newcommand{\qed}{\hbox to5pt{ } \hfill \openbox\bigskip\medskip}

\newcommand{\Fp}{\mathbb F _p}
\newcommand{\Fq}{\mathbb F _q}

\newcommand{\cF}{\mbox{$\cal F$}}
\newcommand{\cG}{\mbox{$\cal G$}}
\newcommand{\cH}{\mbox{$\cal H$}}

\newcommand{\F}{\mathbb F}

\title{A Bollob\'as--type theorem for affine subspaces }
\author{G\'abor Heged\"{u}s\\
{\normalsize \'Obuda University, Antal Bejczy Center for Intelligent Robotics}
}
\begin{document}

\maketitle
\begin{abstract}
Let $W$ denote the $n$-dimensional affine space over the finite field $\Fq$.
%Define
%$$
%\cT:=\{V+\underline{x}:~ \underline{x}\in \Fq,   \textrm{ and }V \textrm{ is a linear subspace of }W\}
%$$
%as the set of affine subspaces of $W$.
We prove here a Bollob\'as--type upper bound in the case of the set of affine subspaces. We give a construction  of a pair of families of affine subspaces, which shows that our result is almost sharp.
\end{abstract}

\section{Introduction}

First we introduce some notation.

In the following let $q=r^{\alpha}$ be a fixed prime power, $n\geq 1$ be a nonnegative integer. Let $W$ denote the $n$-dimensional affine space over the finite field $\Fq$.

%Let 
%$$
%\cS:=\{S\subseteq \Fq^n:~ S \textrm{ is a linear subspace of }W\}
%$$
%denote the set of linear subspaces of the vector space $W$.
 
%Define
%$$
%\cT:=\{S:~ S\textrm{ is an affine subspace of }W\}
%$$
%as the set of affine subspaces of the space $W$.

%Let 
%$$
%\cL:=\{L\in \cS:~ \mbox{dim }L=1\}
%$$
%denote the set of all {\em lines} of the vector space $W$.

B. Bollob\'as proved in \cite{B} the  following famous result.
\begin{thm}
Let $A_1,\ldots,A_m$ and $B_1,\ldots,B_m$ be two families of sets such that 
$A_i\cap B_j=\emptyset$ if and  only if $i=j$. Then
$$
\sum_{i=1}^m \frac{1}{{|A_i|+|B_i| \choose |A_i|}}\leq 1.
$$
In particular if $|A_i|=r$ and $|B_i|=s$ for each $1\leq i\leq m$, then 
$$
m\leq {r+s \choose r}.
$$
\end{thm}

The following strengthening of the uniform version of Bollob\'as's theorem 
%called the {\em Skew Bollob\'as's theorem} 
was proved by L. Lov\'asz in \cite{L1} using tensor product methods. 
\begin{thm} \label{Lovasz}
If $\cF=\{A_1,\ldots,A_m\}$ is an $r$-uniform family and $\cG=\{B_1,\ldots,B_m\}$ is an $s$-uniform family such that
$$
(a)\ A_i\cap B_i=\emptyset\ 
$$
for each $1\leq i\leq m$ and
$$ 
(b)\ A_i\cap B_j\neq\emptyset
$$
if $i<j\ (1\leq i,j\leq m)$, then
$$
m\leq {r+s \choose r}.
$$
\end{thm}

L. Lov\'asz also proved the following generalization of Bollob\'as' theorem for subspaces of a vector space in \cite{L2}:
\begin{thm} \label{subspace}
Let $\F$ be an arbitrary field and $V$ be an $n$-dimensional vector space over the field $\F$. 

Let $U_1,\ldots,U_m$ denote $r$-dimensional subspaces of $V$ and $V_1,\ldots,V_m$ denote $s$-dimensional subspaces of the vector space $V$. Assume that
$$
(a)\ U_i\cap V_i=\{\underline{0}\}
$$
for each $1\leq i\leq m$ and
$$
(b)\ U_i\cap V_j\neq \{\underline{0}\}
$$
whenever $i<j\ (1\leq i,j\leq m)$. Then
$$
m\leq {r+s \choose r}.
$$
\end{thm}

In the following we give an affine version of Theorem \ref{subspace}. 

We say that a pair of families of affine subspaces $(A_i,B_i)_{1\leq i\leq m}$ of $W$ is {\em cross--intersecting} if 
$$
1. \ A_i\cap B_i=\emptyset,
$$
for each $1\leq i\leq m$ and
$$
2. \ A_i\cap B_j\neq \emptyset
$$
whenever $i<j,$\  $(1\leq i,j\leq m)$.

Let $m(n,q)$ denote the maximal size of a cross--intersecting pair of families of affine  subspaces $(A_i,B_i)_{1\leq i\leq n}$.

Our main result is the following modification of Lov\'asz' Theorem \ref{subspace}:

\begin{thm} \label{main}
Let $A_1,\ldots,A_m$ and $B_1,\ldots,B_m$ be  affine subspaces of an $n$-dimensional affine space $W$ over the finite field $\Fq$, where $q\neq 2$. Assume that
$(A_i,B_i)_{1\leq i\leq m}$ is cross--intersecting.
Then 
$$
m\leq q^n+1,
$$
\end{thm}
\begin{remark} Theorem \ref{main} means that
$$
m(n,q)\leq q^n+1.
$$
\end{remark}
\begin{remark}
Our result is a strengthening of Theorem \ref{Lovasz} in the case of affine  hyperplanes. %Namely it follows from  Theorem \ref{Lovasz} that
%if $(A_i,B_i)_{1\leq i\leq m}$  is a cross--intersecting pair of affine translated {\em hyperplanes}, then $m\leq {2q^{n-1} \choose q^{n-1}}$.
\end{remark}

In Section 2 we prove Theorem \ref{main}.
In the proof we use  the polynomial subspace method (see \cite{BF}). 

In Section 3 we give a simple construction, which shows that $m(n,q)\geq \frac{q^n-1}{q-1}$.

Finally in Section 4 we collect  some open problems.

\section{The proof of the main result}

We use the following obvious observation in our proof.

\begin{prop} \label{meet}
The intersection of a family of affine subspaces is either empty or equal to a translate of the intersection of their corresponding vector subspaces.
\end{prop} \qed

Recall that our main result was the following: 

\begin{thm} \label{main22}
Let $A_1,\ldots,A_m$ and $B_1,\ldots,B_m$ be  affine subspaces of an $n$-dimensional affine space $W$ over the finite field $\Fq$, where $q\neq 2$. Assume that
$(A_i,B_i)_{1\leq i\leq m}$ is cross--intersecting.
Then 
$$
m\leq q^n+1,
$$
\end{thm}
{\bf Proof:} 

Let $p$ be an arbitrary, but fixed prime divisor of $q-1$. Since $q\neq 2$, hence $p>1$. We can assign for each subset $F\subseteq \Fq^n$ its characteristic vector $\underline{v_F}\in \{0,1\}^{q^n}\subseteq \Fp^{q^n}$ such that $\underline{v_F}(s)=1$ iff $s\in F$.  Here $\underline{v_F}(s)$ denotes the $s^{th}$ coordinate of the vector $\underline{v_F}$.

Let $1\leq j\leq m$ be fixed. Let $\underline{v_j}=(\underline{v_j}(1),\ldots , \underline{v_j}(q^n))$ denote the characteristic vector of the affine subspace $A_j$ and let $\underline{w_j}=(\underline{w_j}(1),\ldots , \underline{w_j}(q^n))$ denote the  characteristic vector of the affine subspace $B_j$. Here $\underline{v_j}(i)$ denotes the $i^{th}$ coordinate of the vector $\underline{v_j}$.  Similarly, $\underline{w_j}(i)$ denotes the $i^{th}$ coordinate of the vector $\underline{w_j}$.

Consider the polynomials
$$
P_i(x_1,\ldots,x_{q^n}):=1-(\sum_{k=1}^{q^n} \underline{v_j}(k) x_k)\in\Fp[x_1,\ldots,x_{q^n}]
$$
for each $1\leq i\leq m$.

We claim that the polynomials $\{P_i:~ 1\leq i\leq m\}$ are linearly independent functions over $\Fp$. Namely 
$$
P_i( \underline{w_i})=1-\sum_{k=1}^{q^n} \underline{v_i}(k) \underline{w_i}(k)=1-|A_i\cap B_i|=1
$$
 and 
\begin{equation} \label{eq3}
P_i( \underline{w_j})=1-\sum_{k=1}^{q^n} \underline{v_i}(k) \underline{w_j}(k)=1-|A_i\cap B_j|=1-q^t,
\end{equation}
where $t\geq 0$, because  $(A_i,B_i)_{1\leq i\leq m}$ is a cross--intersecting pair of families of affine subspaces and hence we can apply Proposition \ref{meet}. Since 
$$
q\equiv 1 \pmod p,
$$
thus
\begin{equation} \label{eq2}
1-q^t\equiv 0 \pmod p.
\end{equation}
%Therefore the $m\times m$ matrix $P=(P_i(v_j))_{1\leq i,j\leq m}$ is upper triangular over $\Fp$ and in the diagonal we find nonzero elements modulo $p$. Consequently the matrix is nonsingular and it follows from the Triangular Criterion (Proposition 2.8 in \cite{BF}) that 
Consider a linear combination 
$$
\sum_{r=1}^m \lambda_r P_r=0,
$$
where $\lambda_r \in \Fp$. It is easy to prove that $\lambda_r =0$ for each 1$\leq r \leq m$. 
Namely for contradiction, suppose that  there exists  a nontrivial linear relation  
\begin{equation} \label{eq}
\sum_{s=1}^m \lambda_s P_s=0.
\end{equation}
Let $s_0$ be the smallest $s$ such that $\lambda_s\ne 0$. Substitute $ \underline{w_{s_0}}$ for the variable of each side of $(\ref{eq})$. Then by equations (\ref{eq3}) and  (\ref{eq2}), all but the $s_0^{th}$ term vanish, and what remains is
$$
\lambda_{s_0} P_{s_0}(\underline{w_{s_0}})=0.
$$

But  $P_{s_0}(\underline{w_{s_0}})\ne 0$ implies that $\lambda_{s_0}=0$, a contradiction. 
Hence the polynomials $P_1,\ldots,P_m$ are linearly independent functions over $\Fp$.

We infer that the linearly independent polynomials $\{P_1,\ldots,P_m\}$ are in the $\Fp$-space spanned by the monomials
$$
\{x^u\in \Fp[x_1,\ldots,x_{q^n}]:~ \mbox{deg}(x^u)\leq 1\}. 
$$

Clearly
$$
|\{x^u:~ \mbox{deg}(x^u)\leq 1\}|\leq q^n+1,
$$
hence 
$$ 
m\leq q^n+1,
$$
which was to be proved. \qed

\section{A simple construction}

We use in our contruction the following simple proposition.

\begin{prop} \label{eltolt}
Let $F_j$ be arbitrary  affine subspaces for each $1\leq j\leq m$. Let $G_j:=\underline{\alpha_j}+F_j$, where $\underline{\alpha_j}\notin F_j$. Then $F_i\cap G_j\ne \emptyset$ iff $\underline{\alpha_j} \in F_i-F_j$.  
\end{prop}

\proof 

First suppose that $\underline{\alpha_j}\in  F_i-F_j$. Then we can write $\underline{\alpha_j}$ into the form
$$
\underline{\alpha_j}=\underline{f_i}-\underline{f_j},
$$
where $\underline{f_i}\in F_i$ and $\underline{f_j}\in F_j$. Hence $\underline{f_i}=\underline{\alpha_j}+\underline{f_j}\in \underline{\alpha_j}+F_j=G_j$.

On the other hand, suppose that $F_i\cap G_j\neq\emptyset$. Let $\underline{v}\in F_i\cap G_j$, i.e., $\underline{v}\in F_i$ and $\underline{v}\in \underline{\alpha_j}+F_j$. Then there exists $\underline{f_j}\in F_j$ such that $\underline{v}=\underline{\alpha_j}+\underline{f_j}$by definition. Hence $\underline{\alpha_j}=\underline{v}-\underline{f_j}\in F_i-F_j$. \qed

\begin{prop} \label{const}
Let $n\geq 1$ and $q$ be an arbitrary prime power. Then  $m(n,q)\geq \frac{q^n-1}{q-1}$.
\end{prop}

\proof 
Let $m=\frac{q^n-1}{q-1}$. We give a concrete cross--intersecting pair of families  of affine  subspaces $\{A_1,\ldots,A_m\}$ and $\{B_1,\ldots,B_m\}$ of an $n$-dimensional affine space $W$ over the finite field $\Fq$. 
Let 
$$
\cH=\{H_1,\ldots,H_m\}
$$
denote an enumeration of the set of hyperplanes of the vector space $\Fq^n$. It is easy to see that  $m=\frac{q^n-1}{q-1}$. For each $1\leq i\leq m$ we fix a vector $\underline{\beta_i}\in \Fq^n\setminus H_i$.  Define
$$
A_i:=H_i,
$$
and
$$
B_i:=H_i+\underline{\beta_i}.
$$

Clearly $A_i, B_i$ are affine  subspaces of $W$ for each $1\leq i\leq m$. 

Since  $\underline{\beta_i}\notin H_i$ for each $1\leq i\leq m$, hence $A_i\cap B_i=\emptyset$ by the definition of $A_i$ and $B_i$.

On the other hand, since $\underline{\beta_i}\in H_i-H_j=\Fq^n$, hence it follows from Proposition \ref{eltolt} that $A_i\cap B_j\neq \emptyset$ for each $1\leq i<j\leq m$. 
\qed

\section{Open problems}

Here we collect some interesting open problems.

{\em Open problem 1}: What can we say about $m(n,2)$?

\smallskip

{\em Open problem 2}: What is the precise value of  $m(n,q)$, if $q>2$?

\smallskip

%{\em Open problem 3}: What can we say about the projective version of Theorem \ref{main}?

\bigskip
Finally we conjecture the following projective version of Theorem \ref{main}:

\begin{conjecture}
Let $\F$ be an arbitrary field. Let $A_1,\ldots,A_m$ and $B_1,\ldots,B_m$ be projective subspaces of an $n$-dimensional projective space $W$ over the field $\F$. Assume that
$(A_i,B_i)_{1\leq i\leq m}$ is cross--intersecting (i.e. $A_i\cap B_i=\emptyset$ for each $1\leq i\leq m$ and
$A_i\cap B_j\neq \emptyset$ whenever $1\leq i<j\leq m$).
Then 
$$
m\leq 2^{n+1}-2.
$$
\end{conjecture}
      
%We plan to analyze Open Problem 3 in a forthcoming publication.


\begin{thebibliography}{MM}
%\bibitem{AL} W. W. Adams and P. Loustaunau, {\it An Introduction to
%Gr\"obner Bases,} American Mathematical Society, 1994.
\bibitem{BF} L. Babai, P. Frankl, {\em Linear algebra methods in
combinatorics,} September 1992.

\bibitem{B} B. Bollob\'as,  On generalized graphs. Acta Mathematica Hungarica, {\bf 16(3)}, (1965) 447-452.


%\bibitem{CCS} A. M. Cohen, H. Cuypers and H. Sterk (eds.), {\it Some
%apas of Computer Algebra,} Springer-Verlag, Berlin, Heidelberg, 1999.

\bibitem{F} Z. F\"uredi, Geometric solution of an intersection problem for two hypergraphs, {\em European J. of Comb.} {\bf 5} (1984) 133--136.


\bibitem{L1} L. Lov\'asz, Flats in matroids and geometric graphs, in: {\em Combinatorial surveys}, Proc. 6th British Comb. Conf., Egham 1977, Acad. Press, London 1977, 45--86.

\bibitem{L2} L. Lov\'asz, (1979). Topological and algebraic methods in graph theory. In Graph theory and related topics (Proc. Conf., Univ. Waterloo, Waterloo, Ont., 1977)  1-14.

\bibitem{PR} P. Pudl\'ak, V. R\"odl, A combinatorial approach to complexity, {\em Combinatorica} {\bf 12} (1992), 221--226.

 \bibitem{T} Zs. Tuza, Application of Set-Pair Method in Extremal Hypergraph Theory, in ``Extremal problems for Finite Sets'', {\em Bolyai Society Mathematical Studies} {\bf 3}, J\'anos Bolyai Math. Soc., Budapest, 1994, 479--514. 
\end{thebibliography}
\end{document}